\documentclass[11pt]{article}
\usepackage{color}
\usepackage{graphicx}
\usepackage{amsmath,amsthm}
\usepackage{amssymb,mathrsfs}
\usepackage{geometry}
\usepackage[format=hang,font=small,textfont=it]{caption}
\usepackage[nottoc]{tocbibind}
\usepackage{verbatim}
\usepackage{hyperref}
\usepackage{soul}
\soulregister\cite7
\soulregister\eqref7
\hypersetup{colorlinks, linkcolor=blue}

\newtheorem{thm}{Theorem}[section]
\newtheorem{lemma}[thm]{Lemma}
\newtheorem{prop}[thm]{Proposition}
\newtheorem{coro}[thm]{Corollary}

\theoremstyle{definition}
\newtheorem{definition}[thm]{Definition}

\newcommand\ddd{\mathrm{d}}
\newcommand\holder{\mathrm{H\ddot{o}lder}}
\newcommand\ora{\overrightarrow}

\newcommand\bR{\mathbb{R}}

\newcommand\bS{\mathbb{S}}

\def \l {\left}
\def \r {\right}

\begin{document}
	
\renewcommand{\thefootnote}{\fnsymbol {footnote}}
\title{Some Weighted Estimates on Gaussian Measure Spaces}
\footnotetext {{}{2010 \emph{Mathematics Subject Classification}: Primary 42B35; Secondary 42B20, 42B25.}} \footnotetext {{}\emph{Key words and phrases}: Local multilinear fractional integral operators, local Muckenhoupt weights,  rough kernel, Gaussian measure spaces.}
\setcounter{footnote}{0}\author{Boning Di, Qianjun He\footnote{Correspoding author}, Dunyan Yan}
\date{}
\maketitle

\begin{abstract}
In this paper, we obtain the weighted boundedness for the local multi(sub)linear Hardy-Littlewood maximal operators and local multilinear fractional integral operators associated with the local Muckenhoupt weights on Gaussian measure spaces. We deal with these problems by introducing a new pointwise equivalent ``radial'' definitions of these local operators. Moreover using a similar approach, we also get the weighted boundedness for the local fractional maximal operators with rough kernel and local fractional integral operators with rough kernel on Gaussian measure spaces.
\end{abstract}

\section{Introduction}
Weighted norm inequalities arise naturally in harmonic analysis. Furthermore, they have a variety of applications in many fields such as the study of boundary value problem for Laplace's equation on Lipschitz domains, extrapolation theory, vector-valued inequalities, and estimates for certain classes of nonlinear partial differential equations, see for instance, Grafakos~\cite{Grafakos_2014} and Stein~\cite{Stein_1970}. For further information and historical remarks about the one-weight and two-weight estimates, see Muckenhoupt~\cite{Muckenhoupt_1972}, Muckenhoupt and Wheeden~\cite{MW_1974}, Sawyer~\cite{Sawyer_1982}, Sawyer and Wheeden~\cite{SW_1992}, $\mathrm{P\acute{e}rez}$~\cite{Perez_1990, Perez_1994} and Cruz-Uribe et al.~\cite{CMP_2007}.

On the other hand, multivariable calculus provides a strong approach into the study of functions of several variables that goes beyond the narrow perspective of studying a single variable by fixing the other ones. The multilinear fractional integrals have been studied by Grafakos~\cite{Grafakos_1992}, Kenig and Stein~\cite{KS_1999}, Grafakos and Kalton~\cite{GK_2001} and Moen~\cite{Moen_2009}, etc. The multilinear maximal operators have been studied by Lerner et al.~\cite{LOPTT_2009}. On the weighted norm inequalities for operators with rough kernel, recall that Kurtz and Wheeden~\cite{KW_1979} have proved the weighted estimates for fractional integrals with rough kernel under the size condition and the Dini smoothness condition, then Watson~\cite{Watson_1990} and Duoandikoetxea~\cite{Javier_1993} removed the smoothness requirement by using Fourier transform methods. Then in 1998, Ding and Lu~\cite{DL_1998} established the weighted norm inequalities for fractional maximal and singular integral operators with rough kernel. For more about this topic, see Garc\'{\i}a-Cuerva and Rubio~\cite{GR_1985}, Lu et al.~\cite{LDY_2007}.

The purpose of this article is to establish some weighted estimates for the local multilinear operators and the local fractional operators with rough kernel on Gaussian measure spaces. The Gaussian measure space $(\bR^d,|\cdot|,\gamma)$ is the Euclidean space $\bR^d$ endowed with the Euclidean distance $|\cdot|$ and the Gaussian measure $\gamma$, where
\[\ddd\gamma(x):=\pi^{-d/2}e^{-|x|^2} \ddd x.\]
Note that the Gaussian measure is a probability measure, highly concentrated around the origin, with exponential decay at inﬁnity. Thus it isn't a doubling measure, i.e., there is no constant $C>0$, independent of the point $x\in \bR^d$ and real number $r>0$ such that
\[\gamma(B(x,2r))\leq C \gamma(B(x,r))\]
holds for all $x\in\bR^d$ and $r>0$, where $\gamma(B):=\int_B \ddd\gamma(x)$. For more details see~\cite[Appendix 10.3]{GHA_2019}. From this we know that the Gaussian measure space is not a homogeneous type space in the sense of Coifman and Weiss~\cite{CW_1977}. Thus there is few overlap between Gaussian harmonic analysis and harmonic analysis of spaces of homogeneous type. Moreover, the Gaussian measure is trivially a $d$-dimensional measure in $\bR^d$, i.e.,
\[\gamma(B(x,r))\leq Cr^d\]
holds for some positive constant $C$. Hence some results on the $d$-dimensional measure may be useful in Gaussian harmonic analysis, for more about this, see~\cite{GM_2001} and~\cite{LSW_2013}.
On Gaussian measure spaces we study the following admissible balls
\begin{align*}
\mathscr{B}_a:&=\l\{B(x,r): 0<r<a\l(1\wedge \frac{1}{|x|}\r)\r\} \\
&=\l\{B(x,r): 0<r<am(x)\r\},
\end{align*}
where $m(x)=1\wedge \frac{1}{|x|}:=\min\{1,\frac{1}{|x|}\}$ and then $\mathscr{B}_a(x):=\{B: B\in\mathscr{B}_a, x\in B\}$.

Recently the study of Gaussian harmonic analysis has received increasing attention. In 2007, Mauceri-Meda~\cite{MM_2007} studied the $\mathrm{BMO}$ and $H^1$ spaces for the Ornstein-Uhlenbeck operator; while Aimar et al.~\cite{AFS_2007} obtained the weak type $(1,1)$ boundedness for the maximal function and Riesz transform on the Gaussian setting. In 2010, Liu and Yang~\cite{LY_2010} obtained the boundedness for the local fractional integral operators defined by
\[I_{\beta}^{a}(f)(x):=\int_{B(x, am(x))} \frac{f(y)}{[\gamma(B(x,|x-y|))]^{1-\beta}}\ddd\gamma(y)\]
and local fractional maximal operators defined by
\[M_{\beta}^a(f)(x):= \sup_{B\in\mathscr{B}_a(x)} \frac{1}{[\gamma(B)]^{1-\beta}}\int_B |f(y)|\ddd\gamma(y)\]
on the Gaussian Lebesgue spaces. We should point out that they have also obtained some other results which are related to the subject of this paper such as \cite{LY_2008, LY_2010_1, LY_2014, LXYY_2018}.
Then in 2014, Liu et al.~\cite{LSY_2014} established the strong type boundedness for the local fractional integrals and local maximal operators on the Gaussian Morrey-type spaces, etc.; later in 2016, Wang et al.~\cite{WZL_2016} characterized the boundedness for the local maximal operators on the weighted Gaussian Lebesgue spaces by the local $A_{p,a}$ weights defined as
\[[\omega]_{A_{p,a}}=\sup_{B\in\mathscr{B}_{a}}\l(\frac{1}{\gamma(B)}\int_B \omega(x) \ddd\gamma(x)\r)\l(\frac{1}{\gamma(B)}\int_B \omega(x)^{1-p'}\ddd\gamma(x)\r)^{p-1}<\infty,\]
and also obtained the weighted boundedness for the fractional maximal operators on Gaussian measure spaces; in 2020, Lin and Mao~\cite{LM_2020} established the weighted norm inequalities associated with the local Muckenhoupt weights for the local fractional integrals on Gaussian measure spaces.

Based on these results above, we investigate the weighted boundedness for the local multi(sub)linear maximal operators and local multilinear fractional integrals, as well as the local fractional integrals with rough kernel with respect to the local Muckenhoupt weights on the Gaussian Lebesgue spaces. These problems are solved by using the new pointwise equivalent ``radial'' definitions of the local operators on Gaussian measure spaces such as
\begin{align*}
M_{\beta}^a(f)(x)&:= \sup_{B\in\mathscr{B}_a(x)} \frac{1}{[\gamma(B)]^{1-\beta}}\int_B |f(y)|\ddd\gamma(y),\\
&:\sim_{a,d} \sup_{B\in\mathscr{B}_a(x)}\frac{e^{-\beta |x|^2}}{r^{d-d\beta}}\int_{B} |f(y)|\ddd y,
\end{align*}
where $r$ is the radius of $B$. Our main method is transforming the operators' integral measure from Gaussian measure into Lebesgue measure, then these operators possess good translation invariant property and radial property in some sense. See Section \ref{WEFLMFI} below for details.

This paper is organized as follows. Section \ref{DHY_2012_preliminaries} contains some elementary definitions and facts about the Gaussian measure, local Muckenhoupt weights, and some local operators studied by other mathematicians before. The readers who are familiar with this subject may skip directly to Section \ref{DHY_2012_main results}, where all the lemmas and theorems are stated. Finally all the proofs of these results are presented in Section \ref{DHY_2012_proofs}.

We end this section with some notions and notations. Hereafter, we will be working in $\bR^d$ and $d$ will always denote the dimension. We will denote by $C$ or its variants a positive constant independent of the main involved parameters, and use $f\lesssim g$ to denote $f\leq C g$; particularly, if $f\lesssim g\lesssim f$, then we will write $f\sim g$. If necessary, we will denote the dependence of the constants parenthetically: e.g., $C=C(a,d)$ or $C=C_{a,d}$. Similarly, $f\lesssim_{a,d} g$ will denote $f\leq C_{a,d} g$ and $f\sim_{a,d}g$ will denote $C_{a,d} f\leq g\leq C_{a,d} f$. By a weight we will always mean a locally integrable function and non-negative almost everywhere on $\bR^d$ (with respect to the associated measure $\gamma$). For a given weight $\omega$, the weighted Gaussian Lebesgue norm on $\bR^d$ will be denoted by
\[\|f\|_{L^p(\bR^d,\gamma, \omega)}:=\l(\int_{\bR^d} |f(x)|^p\omega(x)\ddd\gamma(x)\r)^{\frac{1}{p}}.\]

\section{Preliminaries}\label{DHY_2012_preliminaries}
\subsection{Gaussian measure and admissible balls}
We have shown that the Gaussian measure is not a doubling measure in $\bR^d$. Nevertheless if we concentrate on the family of admissible balls $\mathscr{B}_a$, then \cite[Proposition 2.1]{MM_2007} points out that the Gaussian measure is doubling if we restrict it to $\mathscr{B}_a$. In other words, there exists a constant $C=C_{a,d}\geq 1$ such that for all $B=B(c_B, r)\in \mathscr{B}_a$ we have
\[\gamma(B(c_B,2r))\leq C \gamma(B(c_B,r)).\]
Hence we say $\gamma$ satisfies the locally doubling condition on $\mathscr{B}_a$. In addition, by \cite[Proposition 2.1]{MM_2007}, the estimate
\begin{equation}\label{DHY_1}
e^{-2a-a^2}\leq e^{|c_B|^2-|x|^2}\leq e^{2a}.
\end{equation}
holds for all $B\in \mathscr{B}_a$ and all $x\in B$. From this fact we know that
\[\gamma(B)\sim_{a,d} e^{-|c_B|^2}r^d\]
holds for all $B\in \mathscr{B}_a$ and moreover, $\gamma$ satisfies the locally reverse doubling condition on $\mathscr{B}_a$ too. Namely, there exists a constant $C=C_{a,d}>1$ such that
\[\gamma(B(c_B,2r))\geq C \gamma(B(c_B,r))\]
holds for all $B=B(c_B, r)\in \mathscr{B}_a$.

\subsection{Local Muckenhoupt weights and Local operators}
Due to the special geometric properties on Gaussian measure spaces, we study the following local $A_{p,a}$ weights and local $A_{p,q,a}$ weights, which are first introduced by Wang et al.\cite{WZL_2016}.
\begin{definition}
Let $a\in(0,\infty)$, $p\in[1,\infty)$. We say a weight $\omega\in A_{p,a}$, if there exists a constant $C\in [1,\infty)$ such that
\[\l(\frac{1}{\gamma(B)}\int_B \omega(x) \ddd\gamma(x)\r)\l(\frac{1}{\gamma(B)}\int_B \omega(x)^{1-p'}\ddd\gamma(x)\r)^{p-1}\leq C\]
holds for all balls $B\in \mathscr{B}_a$; In the case $p=1$,
\[\frac{1}{\gamma(B)}\int_B \omega(x)\ddd\gamma(x)\leq C\; \mathrm{ess}\inf_{x\in B} \omega(x)\]
holds for all balls $B\in \mathscr{B}_a$. Finally, we define the $A_{\infty,a}$ class as $A_{\infty,a}=\cup_{p\geq1}A_{p,a}$.
\end{definition}

\begin{definition}
Let $a>0$ and $1<p,q<\infty$. We say a weight $\omega\in A_{p,q,a}$, if there exists a constant $C>0$ such that 
\[\l(\frac{1}{\gamma(B)}\int_B \omega(x)^q\ddd \gamma(x)\r)^{1/q}\l(\frac{1}{\gamma(B)}\int_B \omega(x)^{-p'}\ddd \gamma(x)\r)^{1/{p'}}\leq C\]
holds for all balls $B\in \mathscr{B}_a$; and say $\omega\in A_{1,q,a}$ $(1<q<\infty)$, if there exists a constant $C>0$ such that
\[\l(\frac{1}{\gamma(B)}\int_B \omega(x)^q\ddd \gamma(x)\r)^{1/q}\l(\mathrm{ess} \sup_B \frac{1}{\omega(x)}\r)\leq C\]
holds for all balls $B\in \mathscr{B}_a$.
\end{definition}
\noindent Note that the class $A_{p,q,a}$ sometimes may be written as $A(p,q,a)$.

It's well known that the Muckenhoupt weights $A_p$ characterize the weighted boundedness for the Hardy-Littlewood maximal operators~\cite{Muckenhoupt_1972}. Similarly, on Gaussian measure spaces~\cite{WZL_2016}, the local Muckenhoupt weights $A_{p,a}$ characterize the weighted boundedness for the local Hardy-Littlewood maximal operators $M_a$ on Gaussian Lebesgue spaces defined by
\[M_a (f)(x):=\sup_{B\in \mathscr{B}_a(x)}\frac{1}{\gamma(B)} \int_B |f(y)| \ddd\gamma(y).\]
Moreover, the local $A_{p,q,a}$ weights give the weighted boundedness for local fractional Hardy-Littlewood maximal operators $M_{\beta}^a$ on Gaussian Lebesgue spaces defined by
\[M_{\beta}^a(f)(x):= \sup_{B\in\mathscr{B}_a(x)} \frac{1}{[\gamma(B)]^{1-\beta}}\int_B |f(y)|\ddd\gamma(y).\]

\noindent Finally, the local Hardy-Littlewood maximal operators $M_{\nu,a}$ with measure $\nu\ddd\gamma$ are defined by
\[M_{\nu,a}(f)(x):=\sup_{B\in\mathscr{B}_a(x)}\frac{1}{\nu(B)}\int_B |f(y)|\nu(y) \ddd\gamma(y),\]
where $\nu(B):=\int_B \nu(x)\ddd\gamma(x)$. These local maximal operators are studied by Lin et al.\cite{LSW_2013} and will be used later.

\section{Main results}\label{DHY_2012_main results}
\subsection{Weighted estimates for local multi(sub)linear maximal operators}
We consider the local multi(sub)linear maximal operators on Gaussian measure spaces first, which are in some way more similar to the classical operators on Lebesgue spaces.
\begin{definition}
Let $f_i\in L_{\mathrm{loc}}^1(\gamma)$, $i=1,2,\ldots,m$  and $a\in(0,\infty)$. The \textit{local multi(sub)linear Hardy-Littlewood maximal operator $\mathscr{M}_a$ on Gaussian measure spaces} is defined by
\[\mathscr{M}_a(f_1,f_2,\ldots,f_m)(x):=\sup_{B\in\mathscr{B}_a(x)} \prod_{j=1}^{m}\frac{\int_{B} |f_j(y)|\ddd \gamma(y)}{\gamma(B)}.\]
\end{definition}
For the multi(sub)linear maximal operators $\mathscr{M}_a$, we investigate the two-weight inequalities and the one-weight inequalities. Our main result about the two-weight boundedness for the local multi(sub)linear maximal operators is the following theorem.
\begin{thm}\label{LOPTT_2009_thm3.3}
Let $1\leq p_j<\infty$ for $j=1,2,\ldots,m$ and $1/p=1/{p_1}+1/{p_2}+\cdots+1/{p_m}$. If $\nu$ and $\omega_j$ be weights satisfying the condition
\begin{equation}\label{LOPTT_2009_thm3.3_0.5}
\sup_{B\in\mathscr{B}_a}\l[\frac{1}{\gamma(B)}\int_B \nu(x) \ddd \gamma(x)\r]^{1/p}\cdot\prod_{j=1}^{m}\l[\frac{1}{\gamma(B)} \int_{B} \omega_j(x)^{1-p'_j}\ddd \gamma(x)\r]^{1/{p'_j}}<\infty.
\end{equation}
In the case $p_j=1$, $\l[\frac{1}{\gamma(B)} \int_{B} \omega_j(x)^{1-p'_j}\ddd \gamma(x)\r]^{1/{p'_j}}$ is understood as $\l(\mathrm{ess} \inf_B \omega_j\r)^{-1}$. In addition, the measure $\nu(x)\ddd\gamma(x)$ satisfies the $a$-local $5$-condition, i.e.
\[\int_{5B}\nu(x)\ddd \gamma(x)\lesssim \int_B \nu(x)\ddd \gamma(x)\]
holds for all balls $B\in \mathscr{B}_a$. Then the following estimate
\[\l\|\mathscr{M}_a(f_1,\ldots,f_m)\r\|_{L^{p,\infty}(\bR^d, \gamma, \nu)}\lesssim \prod_{j=1}^{m}\|f_j\|_{L^{p_j}(\bR^d, \gamma, \omega_j)}\]
holds for all $f_j\in L_{\mathrm{loc}}^1(\gamma)$.
\end{thm}
In order to present the one-weight estimate, we introduce the local multiple weights $A_{\ora{P},a}$ on Gaussian measure spaces as follows.
\begin{definition}
Let $1\leq p_1,p_2,\ldots,p_m<\infty$, $1/p=1/{p_1}+\ldots+1/{p_m}$ and vector $\ora{P}=(p_1,\ldots,p_m)$. Given $\overrightarrow{\omega}=(\omega_1,\omega_2,\ldots,\omega_m)$, we say that $\overrightarrow{\omega}$ satisfies the \textit{$A_{\overrightarrow{P},a}$ condition} if
\[\sup_{B\in\mathscr{B}_a} \l(\frac{1}{\gamma(B)}\int_B \prod_{j=1}^m \omega_j^{p/{p_j}} \ddd\gamma\r)^{1/{p}}\cdot \prod_{j=1}^m \l(\frac{1}{\gamma(B)}\int_B \omega_j^{1-p'_j}\ddd\gamma\r)^{1/{p'_j}}<\infty.\]
When $p_j=1$, the term $\l(\frac{1}{\gamma(B)}\int_B \omega_j^{1-p'_j}\ddd\gamma\r)^{1/{p'_j}}$ is understood as $(\mathrm{ess} \inf_B \omega_j)^{-1}$.
\end{definition}
Based on the definition above we can obtain the following characterization of multilinear $A_{\ora{P},a}$ weights in terms of the linear $A_{p,a}$ weights.
\begin{thm}\label{LOPTT_2009_thm3.6}
Let $\overrightarrow{\omega}=(\omega_1,\omega_2,\ldots,\omega_m)$, $1\leq p_1,p_2,\ldots,p_m<\infty$ and $\nu_{\overrightarrow{\omega}}:=\prod_{j=1}^{m}\omega_j^{p/p_j}$.

Then $\overrightarrow{\omega}\in A_{\overrightarrow{P},a}$ if and only if
\[\l\{\begin{array}{ll}
{\omega_j^{1-p'_j}\in A_{mp'_j, a},} & {j=1,2,\ldots,m}, \\
{\nu_{\overrightarrow{\omega}}\in A_{mp, a},} & {}
\end{array}\right.\]
where the condition $\omega_j^{1-p'}\in A_{mp'_j, a}$ in the case $p_j=1$ is understood as $\omega_j^{1/m}\in A_{1,a}$.
\end{thm}
The proof of this theorem follows from \cite[Theorem 3.6]{LOPTT_2009}, which is in the classical situation. Note that the proof in \cite{LOPTT_2009} mainly uses the $\mathrm{H\ddot{o}lder}$'s inequality that is also valid in Gaussian measure spaces. We omit the detailed proof, because it is too long but essentially the same as \cite[Theorem 3.6]{LOPTT_2009}.

Finally we state the main result about the one-weight boundedness for the local multi(sub)linear maximal operators on Gaussian measure spaces.
\begin{thm}\label{LOPTT_2009_thm3.7}
Let $1<p_j<\infty$ for $j=1,2,\ldots,m$ and $1/p=1/{p_1}+1/{p_2}+\cdots+1/{p_m}$. Then the inequality
\begin{equation}\label{LOPTT_2009_thm3.7_1}
\l\|\mathscr{M}_a(f_1,f_2,\ldots,f_m)\r\|_{L^{p}(\bR^d,\gamma,\nu_{\overrightarrow{\omega}})}\lesssim\prod_{j=1}^{m}\l\|f_j\r\|_{L^{p_j}(\bR^d,\gamma,\omega_j)}
\end{equation}
holds for all $\overrightarrow{f}:=(f_1,f_2,\ldots,f_m)$ if and only if $\overrightarrow{\omega}$ satisfies the $A_{\overrightarrow{P},a}$ condition.
\end{thm}

\subsection{Weighted estimates for local multilinear fractional integrals}\label{WEFLMFI}
Recall that Iida et al.\cite{IKS_2010} point out when the multilinear fractional integrals are defined by
\[I_{\alpha}^{m}(f_1,f_2,\ldots,f_m)(x):=\int_{\bR^d}\prod_{j=1}^m f_j(x-y_j)|(y_1,y_2,\ldots,y_m)|^{\alpha-mn}\ddd y_1\ddd y_2 \ldots \ddd y_m,\]
with some common hypotheses, the weighted boundedness for
\[L^{p_1}(\bR^d, \mu, \omega_1)\times L^{p_2}(\bR^d, \mu, \omega_2)\times\cdots \times L^{p_m}(\bR^d, \mu, \omega_m)\to L^{p}(\bR^d, \mu, \omega^p)\]
is not always true, where $\mu$ denotes the Lebesgue measure. Nevertheless, they obtain the weighted boundedness for multilinear fractional integrals defined by
\[I_{\alpha,\theta}^m (f_1,f_2,\ldots,f_m)(x):=\int_{\bR^d} \prod_{j=1}^m f_j(x-\theta_j y)|y|^{\alpha-n}\ddd y.\]
Inspired by \cite{IKS_2010}, we introduce a new definition for the local fractional integral operators, which is pointwise equivalent to the previous local fractional integrals and is ``centered-radial'' in some sense.

It's easy to see that if $|x-y|<am(x)$, then from the estimate (\ref{DHY_1}) we have
\[e^{-|x|^2}\sim_a e^{-|y|^2};\]
and moreover for every $B=B(x,r)\in \mathscr{B}_a$ we have the estimate
\[\gamma(B)\sim_{a,d} e^{-|x|^2} r^d.\]
From these facts and some calculation we conclude the following
\begin{align*}
\int_{B(x, am(x))} \frac{f(y)}{[\gamma(B(x,|x-y|))]^{1-\beta}}\ddd\gamma(y)
&=\int_{B(x,am(x))} \frac{\pi^{-d/2}f(y)e^{-|y|^2}}{[\gamma(B(x,|x-y|))]^{1-\beta}} \ddd y \\
&\sim_{a,d} \int_{B(x,am(x))} \frac{\pi^{-d/2}f(y)e^{-|x|^2}}{[e^{-|x|^2}|x-y|^d]^{1-\beta}} \ddd y \\
&\sim \pi^{-d/2} e^{-\beta |x|^2}\int_{B(x,am(x))} \frac{f(y)}{|x-y|^{d-d\beta}} \ddd y \\
&\sim\pi^{-d/2} e^{-\beta |x|^2}\int_{B(0,am(x))} \frac{f(x-y)}{|y|^{d-d\beta}} \ddd y
\end{align*}
when $f$ is nonnegative. Since we study the weighted estimates on Gaussian measure spaces, we only need to focus on the nonnegative functions. Now we can give the pointwise equivalent radial definition for the previous local fractional integral operators $I_{\beta}^a$ on nonnegative functions in Gaussian measure spaces.
\begin{definition}\label{DHY_2012_PEDLFI}
Let $a\in(0,\infty)$, $\beta\in(0,1)$ and $f\in L_{\mathrm{c}}^{\infty}(\gamma)$ be nonnegative. The \textit{local fractional integral operator $I_{\beta}^{a}$ on Gaussian measure spaces} is defined by
\begin{align*}
I_{\beta}^{a}(f)(x)&:=\int_{B(x, am(x))} \frac{f(y)}{[\gamma(B(x,|x-y|))]^{1-\beta}}\ddd\gamma(y),\\
&:\sim \pi^{-d/2} e^{-\beta |x|^2}\int_{B(0,am(x))} \frac{f(x-y)}{|y|^{d-d\beta}} \ddd y,
\end{align*}
where $L_{\mathrm{c}}^{\infty}(\gamma)$ denotes the set of all functions in $L^{\infty}$ with compact support on $\bR^d$.
\end{definition}

Based on this new definition \ref{DHY_2012_PEDLFI} and \cite{IKS_2010}, on Gaussian measure spaces, we study the following local multilinear fractional integral operators.
\begin{definition}
Let $a\in(0,\infty)$ and $\beta\in (0,1)$. If $\theta_j\neq0$ for every $1\leq j\leq m$ and $\theta_j$ all distinct, we define the \textit{local m-linear fractional integral $I_{\beta,\theta,m}^a$ on Gaussian measure spaces} by
\[I_{\beta,\theta,m}^a (f_1,f_2,\ldots,f_m)(x):=\pi^{-d/2}e^{-\beta|x|^2}\int_{B(0,am(x))} \frac{\prod_{j=1}^m f_j(x-\theta_j y)}{|y|^{d-d\beta}}\ddd y,\]
where $f_j\in L_{\mathrm{c}}^{\infty}(\gamma)$ for all $j=1,\ldots,m$.
\end{definition}
Then we can obtain the following pointwise inequality and futhermore get the weighted estimates for $I_{\beta,\theta,m}^a$ by following the proof in \cite{IKS_2010}.
\begin{lemma}\label{IKS_2010_lemma2.2}
Let $a\in(0,\infty)$, $\beta\in(0,1)$ and $p_j>1 $ with $ j=1,2,\ldots,m$. If $1/s=1/{p_1}+1/{p_2}+\cdots+1/{p_m}$, then we have
\[\l|I_{\beta,\theta,m}^a (f_1,f_2,\ldots,f_m)(x)\r|\lesssim \prod_{j=1}^m \l[I_{\beta}^{\theta_j a}\l(|f_j|^{\frac{p_j}{s}}\r)(x)\r]^{\frac{s}{p_j}}, \quad \forall x\in \bR^d.\]
\end{lemma}
\begin{thm}\label{IKS_2010_thm1}
Let $1/s=1/{p_1}+1/{p_2}+\cdots+1/{p_m}$, $p_j>1$ for $j=1,2,\ldots,m$,  $\beta\in(0,1)$ and $1<s<1/{\beta}$. If $1/p=1/s-\beta$ and $\omega_j(x)^{\frac{p_j}{s}}\in A_{s,p,a}$ for every $j$, then we have
\[\l\|I_{\beta,\theta,m}^{a}(f_1,f_2,\ldots,f_m)\r\|_{L^{p}(\bR^d, \gamma, \omega^p)}\lesssim \prod_{j=1}^{m} \|f_j\|_{L^p(\bR^d,\gamma, \omega_j^{p_j})},\]
where $\omega(x):=\prod_{j=1}^m \omega_j(x)$.
\end{thm}

\subsection{Weighted estimates for local fractional integral operators with rough kernel}
By using a similar method as we give the new radial version of $I_{\beta}^a$ in Definition \ref{DHY_2012_PEDLFI}, we obtain the weighted estimates for local fractional maximal operators with rough kernel on Gaussian measure spaces. Then just like the classical situation in \cite{DL_1998}, we obtain the weighted estimates for local fractional integral operators with rough kernel on Gaussian measure spaces.
\begin{definition}
 Suppose $f\in L_{\mathrm{loc}}^1(\gamma)$, $a\in(0,\infty)$, $\beta\in(0,1)$ and $\Omega(x')\in L^s(\bS^{d-1})$, where $x'=x/|x|$. The \textit{local fractional maximal operator with rough kernel $M_{\Omega, \beta}^a$ on Gaussian measure spaces} is defined by
\[M_{\Omega, \beta}^a(f)(x):= \sup_{B\in\mathscr{B}_a(x)} \frac{1}{[\gamma(B)]^{1-\beta}}\int_B |\Omega(x-y)||f(y)| \ddd\gamma(y).\]
\end{definition}
In order to prove the desired weighted boundedness for $M_{\Omega,\beta}^a$, we first introduce the following basic two results.
\begin{lemma}\label{DL_1998_lemma1}
Let $0<\beta<1$, $1<s'$ and $1<p/{s'}<1/{\beta}$. If $\frac{1}{q/{s'}}=\frac{1}{p/{s'}}-\beta$ and $\omega^{s'}\in A(p/{s'},q/{s'},a)$, then
\[\l\| N_{\beta,s'}^a f\r\|_{L^q(\bR^d,\gamma, \omega^q)}\lesssim \|f\|_{L^p(\bR^d,\gamma,\omega^p)},\]
where $N_{\beta, s'}^a f(x)$ is the fractional maximal operator of order $s'$ with Gaussian measure defined by
\[N_{\beta, s'}^a f(x):= \sup_{B\in\mathscr{B}_a(x)}\l(\frac{1}{\gamma(B)^{1-\beta}}\int_B |f(y)|^{s'} \ddd\gamma(y)\r)^{1/{s'}}.\]
\end{lemma}
\begin{prop}\label{DL_1998_prop1}
Let $0<\beta<1$, $s'<p<1/{\beta}$ and $1/q=1/p-\beta$. If $\Omega\in L^{s}(\bS^{d-1})$ with homogeneity of degree $0$ and $\omega^{s'}\in A(p/{s'},q/{s'},a)$, then
\[\l\|M_{\Omega,\beta}^a f\r\|_{L^q (\bR^d,\gamma, \omega^q)}\lesssim \|f\|_{L^p(\bR^d, \gamma,\omega^p)}.\]
\end{prop}
Based on the proposition \ref{DL_1998_prop1} above, now we give the definition of the local fractional integral operators with rough kernel on Gaussian measure spaces.
\begin{definition}
Let $a\in(0,\infty)$, $\beta\in(0,1)$, $\Omega\in L^s(\bS^{d-1})$ and $f\in L_{\mathrm{c}}^{\infty}(\gamma)$. The \textit{local fractional integral operator with rough kernel $I_{\Omega, \beta}^a$ on Gaussian measure spaces} is defined by
\[I_{\Omega,\beta}^a (f)(x):=\int_{B(x,am(x))} \frac{\Omega(x-y)f(y)}{[\gamma(B(x,|x-y|))]^{1-\beta}}\ddd\gamma(y),\]
where $L_{\mathrm{c}}^{\infty}(\gamma)$ denotes the set of all functions in $L^{\infty}$ with compact support on $\bR^d$.
\end{definition}
In order to get the desired weighted boundedness for local fractional integral operators with rough kernel $I_{\Omega,\beta}^a$ on Gaussian measure spaces, we need to prove the following results one by one.
\begin{lemma}\label{DL_1998_lemma2}
Let $0<\beta<1$, $1<p<1/{\beta}$ and $1/q=1/p-\beta$. If $\omega\in A(p,q,a)$, then there exists an $\varepsilon>0$ such that

(i) $\varepsilon<\beta<\beta+\varepsilon<1$;

(ii) $1/p>\beta+\varepsilon, 1/q<1-\varepsilon$;

(iii) $\omega\in A(p,q_{\varepsilon},a)$ and $\omega\in A(p,\tilde{q}_{\varepsilon},a)$,\\
where $1/{q_{\varepsilon}}=1/p-\beta-\varepsilon$ and $1/{\tilde{q}_{\varepsilon}}=1/p-\beta+\varepsilon$.
\end{lemma}
\begin{lemma}\label{DL_1998_lemma3}
Let $0<\beta<1$, $1\leq s'<p<1/{\beta}$ and $1/q=1/p-\beta$. If $\omega^{s'}\in A(p/{s'},q/{s'},a)$, then there exists an $\varepsilon>0$ such that

(i) $\varepsilon<\beta<\beta+\varepsilon<1$;

(ii) $1/p>\beta+\varepsilon, 1/q<1-\varepsilon$;

(iii) $\omega^{s'}\in A(p/{s'},q_{\varepsilon}/{s'}, a)$ and $\omega^{s'}\in A(p/{s'},\tilde{q}_{\varepsilon}/{s'},a)$,\\
where $1/{q_{\varepsilon}}=1/p-\beta-\varepsilon$ and $1/{\tilde{q}_{\varepsilon}}=1/p-\beta+\varepsilon$.
\end{lemma}
These two lemmas are almost the same as \cite[Lemma 2 and Lemma 3]{DL_1998}. Recall that the proofs in \cite{DL_1998} only used some basic properties of $A_p$ and $A_{p,q}$, which are also valid for local Muckenhoupt weights $A_{p,a}$ and $A_{p,q,a}$ according to \cite{WZL_2016}. Hence we omit these two proofs here.
\begin{lemma}\label{LM_2020_lemma3.4}
Let $a\in (0,\infty)$, $0\leq \beta_1<\beta<\beta_2\leq 1$ and $f\in L_{\mathrm{loc}}^1 (\gamma)$. Then
\[|I_{\Omega,\beta}^{a} (f)(x)|\lesssim [M_{\Omega,\beta_1}^{2a} (f)(x)]^{\frac{\beta_2-\beta}{\beta_2-\beta_1}}[M_{\Omega,\beta_2}^{2a} (f)(x)]^{\frac{\beta-\beta_1}{\beta_2-\beta_1}}, \quad \forall x\in \bR^d.\]
\end{lemma}
From this important pointwise inequality, we can directly get the following corollary
\begin{coro}\label{DL_1998_lemma4}
For any $\varepsilon>0$ with $0<\beta-\varepsilon<\beta+\varepsilon<1$, we have
\[|I_{\Omega,\beta}^a (f)(x)|\lesssim_{\varepsilon,\beta,d} [M_{\Omega,\beta+\varepsilon}^{2a} (f)(x)]^{1/2}[M_{\Omega,\beta-\varepsilon}^{2a} (f)(x)]^{1/2}, \quad \forall x\in \bR^d.\]
\end{coro}
Finally, by using the results above, we can obtain the following desired theorem.
\begin{thm}\label{DL_1998_thm1}
Let $0<\beta<1$, $s'<p<1/\beta$ and $1/q=1/p-\beta$. If
\[\Omega\in L^s(\bS^{d-1}) \quad\text{and}\quad \omega(x)^{s'}\in A\l(\frac{p}{s'}, \frac{q}{s'}, a\r),\]
then
\[\|I_{\Omega,\beta}^{a} f\|_{L^q(\bR^d, \gamma, \omega^q)}\lesssim \|f\|_{L^p(\bR^d, \gamma, \omega^p)}.\]
\end{thm}

\section{Proofs of the main results}\label{DHY_2012_proofs}
\begin{proof}[\textbf{Proof of Theorem \ref{LOPTT_2009_thm3.3}}]
The proof is similar to that in \cite[Theorem 3.3]{LOPTT_2009}. We only present the case where $p_j>1$ for all $j=1,2,\ldots,m$ as the case when some $p_j=1$ is a minor modification of the linear case. Without loss of generality, we can assume $f_j\geq0$ for all  $j$. Let
\[E_{\lambda}=\{x\in\bR^d: |\mathscr{M}_a(f_1,f_2,\ldots,f_m)(x)|>\lambda\}.\]
We need to prove that the following estimate
\begin{equation}\label{LOPTT_2009_thm3.3_1}
\l[\int_{E_{\lambda}} \nu(x)\ddd \gamma(x)\r]^{1/p}\lesssim \lambda^{-1}\prod_{j=1}^{m}\|f_j\|_{L^{p_j}(\bR^d, \nu, \omega_j)} .
\end{equation}
holds for every $\lambda>0$ and every $f_j\in L_{\mathrm{loc}}^1(\gamma)$.

First let's consider the necessary condition for this estimate (\ref{LOPTT_2009_thm3.3_1}). This part will be used to prove Theorem \ref{LOPTT_2009_thm3.7} later, although it may be useless for the proof of Theorem \ref{LOPTT_2009_thm3.3} here. Assuming that the inequality (\ref{LOPTT_2009_thm3.3_1}) is right, for any ball $B\in \mathscr{B}_a$, from the definition we know that
\[\mathscr{M}_a(f_1\chi_B,f_2\chi_B,\ldots,f_m\chi_B)(x)\geq \prod_{j=1}^{m}\frac{1}{\gamma(B)} \int_B f_j(y)\ddd \gamma(y)=:\lambda_0\]
holds for all $x\in B$. Thereby if $\lambda<\lambda_0$, then it's obvious that
\[B\subset \{x\in\bR^d: |\mathscr{M}_a(f_1\chi_B,f_2\chi_B,\ldots,f_m\chi_B)(x)|>\lambda\}=: E_{\lambda, B}.\]
Thus we conclude the following
\[\l(\int_B \nu(x)\ddd\gamma(x)\r)^{1/p}\leq\l(\int_{E_{\lambda, B}} \nu(x)\ddd\gamma(x)\r)^{1/p}\lesssim \lambda^{-1} \prod_{j=1}^{m}\l\|f_j\r\|_{L^{p_j}(B, \gamma, \omega_j)}.\]
Since this inequality holds for all $\lambda<\lambda_0$, letting $\lambda\to\lambda_0$ it follows that
\begin{equation}\label{LOPTT_2009_thm3.3_2}
\prod_{j=1}^{m}\frac{1}{\gamma(B)} \int_B f_j(y)\ddd \gamma(y) \l(\int_B \nu(x)\ddd\gamma(x)\r)^{1/p} \lesssim \prod_{j=1}^{m}\l\|f_j\r\|_{L^{p_j}(B, \gamma, \omega_j)}
\end{equation}
holds for all balls $B\in\mathscr{B}_a$. 

Here we have shown that $(\ref{LOPTT_2009_thm3.3_1})\Rightarrow(\ref{LOPTT_2009_thm3.3_2})$ which will be used in the proof of Theorem \ref{LOPTT_2009_thm3.7} below. Now we turn to the proof of Theorem \ref{LOPTT_2009_thm3.3}. We shall prove this theorem by following the frame:
\[\text{ $a$-local $5$-condition} + (\ref{LOPTT_2009_thm3.3_2})\Rightarrow (\ref{LOPTT_2009_thm3.3_1}) \quad \text{and} \quad (\ref{LOPTT_2009_thm3.3_2}) \Leftrightarrow (\ref{LOPTT_2009_thm3.3_0.5}).\]
Suppose the $a$-local $5$-condition and the estimate (\ref{LOPTT_2009_thm3.3_2}) are right. From (\ref{LOPTT_2009_thm3.3_2}) we can obtain
\[\prod_{j=1}^m \frac{1}{\gamma(B)}\int_B f_j(y)\ddd\gamma(y)\lesssim \frac{\prod_{j=1}^m \l\|f_j\r\|_{L^{p_j}(B, \gamma, \omega_j)}}{\l(\int_B\nu(x)\ddd\gamma(x)\r)^{1/p}}=\prod_{j=1}^m \l[\frac{\int_B f_j(x)^{p_j}\omega_j(x)\ddd\gamma(x)}{\int_B \nu(x)\ddd\gamma(x)}\r]^{1/{p_j}}.\]
For $x\in\bR^d$, by taking supremum over all balls $B\in\mathscr{B}_a(x)$ we get the following
\[\mathscr{M}_a(f_1,f_2,\ldots,f_m)(x)\lesssim\prod_{j=1}^{m}\l[M_{\nu,a}\l(f_j^{p_j}\omega_j/\nu\r)(x)\r]^{1/{p_j}}.\]
This estimate and the $\mathrm{H\ddot{o}lder}$'s inequality for weak spaces (see \cite[Exercise 1.1.15]{Grafakos_2014}) yield that 
\begin{align*}
\l\|\mathscr{M}_a(f_1,f_2,\ldots,f_m)\r\|_{L^{p,\infty}(\bR^d,\gamma,\nu)}
&\lesssim \l\|\prod_{j=1}^{m}\l[M_{\nu,a}\l(f_j^{p_j}\omega_j/\nu\r)\r]^{1/{p_j}}\r\|_{L^{p,\infty}(\bR^d,\gamma,\nu)} \\
&\leq \prod_{j=1}^{m}\l\|\l[M_{\nu,a}\l(f_j^{p_j}\omega_j/\nu\r)\r]^{1/{p_j}}\r\|_{L^{p_j,\infty}(\bR^d,\gamma,\nu)} \\
&=\prod_{j=1}^{m}\l\| M_{\nu,a}\l(f_j^{p_j}\omega_j/\nu\r)\r\|_{L^{1,\infty}(\bR^d,\gamma,\nu)}^{1/{p_j}}.
\end{align*}
Then by using \cite[Proposition 2.2]{LSW_2013}, we can prove the desired inequality (\ref{LOPTT_2009_thm3.3_1}) as follows
\begin{align*}
\l\|\mathscr{M}_a(f_1,f_2,\ldots,f_m)\r\|_{L^{p,\infty}(\bR^d,\gamma,\nu)}
&\lesssim\prod_{j=1}^{m}\l\| M_{\nu,a}\l(f_j^{p_j}\omega_j/\nu\r)\r\|_{L^{1,\infty}(\bR^d,\gamma,\nu)}^{1/{p_j}} \\
&\lesssim\prod_{j=1}^{m} \l\|f_j^{p_j}\omega_j/\nu\r\|_{L^1(\bR^d,\gamma,\nu)}^{1/p_j} \\
&=\prod_{j=1}^m \|f_j\|_{L^{p_j}(\bR^d,\gamma,\omega_j)}
\end{align*}

Finally, we show that the condition (\ref{LOPTT_2009_thm3.3_0.5}) is equivalent to the inequality (\ref{LOPTT_2009_thm3.3_2}) and therefore we complete the proof. First if (\ref{LOPTT_2009_thm3.3_2}) is right, then by setting $f_j(x)=\omega(x)^{-\frac{1}{{p_j-1}}}$ we get (\ref{LOPTT_2009_thm3.3_0.5}) immediately. On the other hand if (\ref{LOPTT_2009_thm3.3_0.5}) is right, by the $\mathrm{H\ddot{o}lder}$'s inequality we can get the following 
\begin{align*}
&\prod_{j=1}^{m}\frac{1}{\gamma(B)} \int_B f_j(y)\ddd \gamma(y)\cdot \l(\int_B \nu(y)\ddd\gamma(y)\r)^{1/p} \\
&\lesssim \prod_{j=1}^{m}\frac{\int_B f_j(y)\ddd \gamma(y)}{\gamma(B)} \cdot \gamma(B)^{1/p}\prod_{j=1}^{m} \l(\frac{1}{\gamma(B)} \int_B \omega_j(y)^{1-p'_j} \ddd\gamma(y)\r)^{\frac{-1}{p'_j}} \\
&=\prod_{j=1}^{m}\int_B f_j(y)\ddd \gamma(y) \cdot \prod_{j=1}^{m} \l(\int_B \omega_j(y)^{1-p'_j} \ddd\gamma(y)\r)^{\frac{-1}{p'_j}} \\
&\leq \prod_{j=1}^m \l[\l(\int_B f_j(y)^{p_j}w_j(y)\ddd\gamma(y)\r)^{\frac{1}{p_j}}\l(\int_B \omega_j(y)^{\frac{-p'_j}{p_j}}\ddd\gamma(y)\r)^{\frac{1}{p'_j}}\r]\cdot \prod_{j=1}^{m} \l(\int_B \omega_j(y)^{1-p'_j} \ddd\gamma(y)\r)^{\frac{-1}{p'_j}} \\
&=\prod_{j=1}^m \l(\int_B f_j(y)^{p_j}\omega_j(y)\ddd\gamma(y)\r)^{1/{p_j}}.
\end{align*}
Thus we have proved the equivalence between the condition (\ref{LOPTT_2009_thm3.3_0.5}) and the inequality (\ref{LOPTT_2009_thm3.3_2}), and the proof is finished.
\end{proof}

\begin{proof}[\textbf{Proof of Theorem \ref{LOPTT_2009_thm3.7}}]
We adapt some ideas from \cite{LOPTT_2009}. Let's consider the necessity first. From the proof of Theorem \ref{LOPTT_2009_thm3.3} we know that the inequality (\ref{LOPTT_2009_thm3.7_1}) leads to the condition (\ref{LOPTT_2009_thm3.3_0.5}), then by applying this condition (\ref{LOPTT_2009_thm3.3_0.5}) we get the necessity immediately. Thereby we only need to prove the sufficiency with $f_j\geq 0$ for all $j=1,\ldots, m$.

If $\ora{\omega}\in A_{\ora{P},a}$, from Theorem \ref{LOPTT_2009_thm3.6} and \cite[Proposition 2.2]{WZL_2016} we obtain that each $w_j^{-1/{(p_j-1)}}$ satisfies the reverse $\holder$ inequality, i.e., there exist $r_j>1$ and $C>0$ such that 
\begin{equation}\label{LOPTT_2009_thm3.7_1.5}
\l(\frac{1}{\gamma(B)}\int_B \omega_j^{-\frac{r}{p_j-1}}\ddd\gamma\r)^{1/r}\leq C\frac{1}{\gamma(B)}\int_B \omega_j^{-\frac{1}{p_j-1}}\ddd\gamma.
\end{equation}
holds for every $1\leq r\leq r_j$ and every $B\in \mathscr{B}_a$. We claim the following pointwise inequality:
\begin{equation}\label{LOPTT_2009_thm3.7_2}
\mathscr{M}_a(\ora{f})(x)\lesssim \prod_{j=1}^{m} M_{\nu_{\ora{\omega}},a}\l[\l(f_j^{p_j}\omega_j/\nu_{\ora{\omega}}\r)^q\r](x)^{1/{qp_j}}.
\end{equation}
By Theorem \ref{LOPTT_2009_thm3.6} we get $\nu_{\ora{\omega}}\in A_{mp, a}$, then from \cite[Lemma 3.1]{WZL_2016} we know that $\nu\ddd\gamma$ is local doubling on $\mathscr{B}_a$. Hence \cite[Proposition 2.2]{LSW_2013} and the Marcinkiewicz interpolation theorem show the boundedness
\[M_{\nu_{\ora{\omega}},a}:
\l\{\begin{array}{ll}
L^p(\omega\ddd\gamma) \to L^p(\omega\ddd\gamma) & p\in(1,\infty) \\
L^1(\omega\ddd\gamma) \to L^{1,\infty}(\omega\ddd\gamma) & p=1
\end{array}\right. .\]
Therefore by the claim (\ref{LOPTT_2009_thm3.7_2}) and $\holder$'s inequality we have
\begin{align*}
\l(\int_{\bR^d} \l|\mathscr{M}_a(\ora{f})(x)\r|^p\nu_{\ora{\omega}}(x)\ddd\gamma(x)\r)^{1/p}
&\lesssim\l\|\prod_{j=1}^{m} M_{\nu_{\ora{\omega}},a}\l[\l(f_j^{p_j}\omega_j/\nu_{\ora{\omega}}\r)^q\r]^{1/{qp_j}}\r\|_{L^p(\bR^d,\gamma,\nu_{\ora{\omega}})} \\
&\leq \prod_{j=1}^m \l\|M_{\nu_{\ora{\omega}},a}\l[\l(f_j^{p_j}\omega_j/\nu_{\ora{\omega}}\r)^q\r]^{1/{qp_j}}\r\|_{L^{p_j}(\bR^d,\gamma,\nu_{\ora{\omega}})}.
\end{align*}
Then using the boundedness for $M_{\nu_{\ora{\omega}},a}$ we can obtain the desired conclusion
\begin{align*}
\l(\int_{\bR^d} \l|\mathscr{M}_a(\ora{f})(x)\r|^p\nu_{\ora{\omega}}(x)\ddd\gamma(x)\r)^{1/p}
&\lesssim \prod_{j=1}^m \l\|M_{\nu_{\ora{\omega}},a}\l[\l(f_j^{p_j}\omega_j/\nu_{\ora{\omega}}\r)^q\r]^{1/{qp_j}}\r\|_{L^{p_j}(\bR^d,\gamma,\nu_{\ora{\omega}})} \\
&\lesssim \prod_{j=1}^m \l\|\l(f_j^{p_j}\omega_j/\nu_{\ora{\omega}}\r)^{1/{p_j}}\r\|_{L^{p_j}(\bR^d,\gamma,\nu_{\ora{\omega}})} \\
&=\prod_{j=1}^{m} \l\|f_j\r\|_{L^{p_j}(\bR^d,\gamma,\omega_j)}.
\end{align*}
It remains to prove the claim (\ref{LOPTT_2009_thm3.7_2}) now. Notice that in the classical situation \cite[Theorem 3.7]{LOPTT_2009}, this claim can be proved by using the $\mathrm{H\ddot{o}lder}$'s inequality. We omit the proof here since the $\mathrm{H\ddot{o}lder}$'s inequality also works on Gaussian measure spaces.
\end{proof}

\begin{proof}[\textbf{Proof of Lemma \ref{IKS_2010_lemma2.2}}]
By using the $\holder$'s inequality and changes of variable formula, we can prove the desired conclusion as follows
\begin{align*}
\l|I_{\beta,\theta,m}^a (f_1,f_2,\ldots,f_m)(x)\r| &\leq \pi^{-d/2}e^{-\beta|x|^2} \int_{B(0,am(x))} \frac{\prod_{j=1}^m |f_j(x-\theta_j y)|}{|y|^{d-d\beta}}\ddd y \\
&\leq \pi^{-d/2}e^{-\beta|x|^2} \prod_{j=1}^m \l(\int_{B(0, am(x))} |f_j(x-\theta_j y)|^{\frac{p_j}{s}}\frac{1}{|y|^{d-d\beta}}\ddd y\r)^{\frac{s}{p_j}} \\
&=\pi^{-d/2}e^{-\beta|x|^2} \prod_{j=1}^m \l(\int_{B(0, \theta_j am(x))} |f_j(x-y)|^{\frac{p_j}{s}}\frac{1}{\l|y/{\theta_j}\r|^{d-d\beta}} \theta_j^{-d}\ddd y\r)^{\frac{s}{p_j}} \\
&=\prod_{j=1}^m \theta_j^{\frac{-d\beta s}{p_j}} \pi^{-d/2}e^{-\beta|x|^2} \prod_{j=1}^m \l(\int_{B(0, \theta_j am(x))} |f_j(x-y)|^{\frac{p_j}{s}}\frac{1}{\l|y\r|^{d-d\beta}} \ddd y\r)^{\frac{s}{p_j}} \\
&=\prod_{j=1}^m \theta_j^{\frac{-d\beta s}{p_j}} \prod_{j=1}^m \l[I_{\beta}^{\theta_j a}\l(|f_j|^{\frac{p_j}{s}}\r)(x)\r]^{\frac{s}{p_j}}.
\end{align*}
\end{proof}

\begin{proof}[\textbf{Proof of Theorem \ref{IKS_2010_thm1}}]
First by using Lemma \ref{IKS_2010_lemma2.2} and $\holder$'s inequality we conclude that
\begin{align*}
\l(\int_{\bR^d} \l|I_{\beta,\theta,m}^{a}(f_1,f_2,\ldots,f_m)(x)\omega(x)\r|^p\ddd\gamma(x)\r)^{1/p} 
&\lesssim\l(\int_{\bR^d} \l|\prod_{j=1}^m \l[I_{\beta}^{\theta_j a}\l(|f_j|^{\frac{p_j}{s}}\r)(x)\r]^{\frac{s}{p_j}}\omega(x)\r|^p\ddd\gamma(x)\r)^{1/p} \\
&=\l(\int_{\bR^d} \l|\prod_{j=1}^m \l[I_{\beta}^{\theta_j a}\l(|f_j|^{\frac{p_j}{s}}\r)(x)\omega_j(x)^{\frac{p_j}{s}}\r]^{\frac{s}{p_j}}\r|^p\ddd\gamma(x)\r)^{1/p} \\
&\leq \l[\prod_{j=1}^m \l(\int_{\bR^d} \l[I_{\beta}^{\theta_j a}\l(|f_j|^{\frac{p_j}{s}}\r)(x)\omega_j(x)^{\frac{p_j}{s}}\r]^p \ddd\gamma(x)\r)^{\frac{s}{p_j}}\r]^{\frac{1}{p}} \\
&=\prod_{j=1}^m \l\|I_{\beta}^{\theta_j a}\l(|f_j|^{\frac{p_j}{s}}\r)\r\|_{L^p(\bR^d, \gamma, \omega_j^{\frac{pp_j}{s}})}.
\end{align*}
Then by using \cite[Corollary 2.9]{LM_2020} and \cite[Theorem 3.8]{LM_2020} we can obtain the following estimate
\begin{align*}
&\l(\int_{\bR^d} \l|I_{\beta,\theta,m}^{a}(f_1,f_2,\ldots,f_m)(x)\omega(x)\r|^p\ddd\gamma(x)\r)^{1/p} \\
&\lesssim\prod_{j=1}^m \l\|I_{\beta}^{\theta_j a}\l(|f_j|^{\frac{p_j}{s}}\r)\r\|_{L^p(\bR^d, \gamma, \omega_j^{\frac{pp_j}{s}})} \\
&\lesssim\prod_{j=1}^m \l\|f_j^{\frac{p_j}{s}}\r\|^{s/{p_j}}_{L^s(\bR^d, \gamma, \omega_j^{\frac{sp_j}{s}})} \\
&=\prod_{j=1}^m \|f_j\|_{L^{p_j}(\bR^d,\gamma,\omega_j^{p_j})}.
\end{align*}
\end{proof}

\begin{proof}[\textbf{Proof of Lemma \ref{DL_1998_lemma1}}]
From the definitions of $N_{\beta, s'}^a$ and $M_{\beta}^a$ we know that $N_{\beta, s'}^a f(x)=\l[M_{\beta}^a (|f|^{s'})(x)\r]^{1/{s'}}$. Then we conclude that
\begin{align*}
\l\| N_{\beta,s'}^a f\r\|_{L^q(\bR^d,\gamma, \omega^q)}
&=\l(\int_{\bR^d} \l[M_{\beta}^a (|f|^{s'})(x)\r]^{q/{s'}}\omega(x)^q \ddd\gamma(x)\r)^{1/q} \\
&=\l[\l(\int_{\bR^d} \l[M_{\beta}^a (|f|^{s'})(x)\nu(x)\r]^{q/{s'}} \ddd\gamma(x)\r)^{s'/q}\r]^{1/{s'}},
\end{align*}
where $\nu(x)=\omega(x)^{s'}\in A(p/{s'},q/{s'},a)$. By \cite[Theorem 5.3]{WZL_2016}, we have
\begin{align*}
\l(\int_{\bR^d} \l[M_{\beta}^a (|f|^{s'})(x)\nu(x)\r]^{q/{s'}} \ddd\gamma(x)\r)^{s'/q}
&\lesssim \l(\int_{\bR^d} \l[|f(x)|^{s'}\nu(x)\r]^{p/{s'}}\ddd\gamma(x)\r)^{s'/p} \\
&=\l(\int_{\bR^d} |f(x)|^{p}\nu(x)^{p/{s'}}\ddd\gamma(x)\r)^{s'/p}.
\end{align*}
Hence we can get the following desired result
\[\l\| N_{\beta,s'}^a f\r\|_{L^q(\bR^d,\gamma, \omega^q)}\lesssim \l(\int_{\bR^d} |f(x)|^{p}\nu(x)^{p/{s'}}\ddd\gamma(x)\r)^{1/p} =\|f\|_{L^p(\bR^d,\gamma,\omega^p)}.\]
\end{proof}

\begin{proof}[\textbf{Proof of Proposition \ref{DL_1998_prop1}}]
By $\holder$'s inequality we obtain that
\begin{equation}\label{DL_1998_prop1_1}
\l|M_{\Omega, \beta}^a(f)(x)\r|\leq \sup_{B\in\mathscr{B}_a(x)} \frac{1}{\gamma(B)^{1-\beta}}\l(\int_B |\Omega(x-y)|^s \ddd\gamma(y)\r)^{1/s}\l(\int_B |f(y)|^{s'} \ddd\gamma(y)\r)^{1/{s'}}.
\end{equation}
Since $B=B(c_B,r)\in\mathscr{B}_a(x)$, we have $|x-y|<2am(x)$ and  $|x-c_B|<am(x)$. Then it follows that
\[e^{-|y|^2}\sim_a e^{-|x|^2}\sim_a e^{-|c_B|^2}.\]
Therefore we conclude the following relation
\begin{align*}
\l(\int_{B(c_B, r)} |\Omega(x-y)|^s e^{-|y|^2}\pi^{-\frac{d}{2}}\ddd y\r)^{\frac{1}{s}}
&\sim_a \pi^{-\frac{d}{2s}} e^{-\frac{|c_B|^2}{s}} \l(\int_{B(c_B, r)} |\Omega(x-y)|^s \ddd y\r)^{\frac{1}{s}} \\
&\sim\pi^{-\frac{d}{2s}} e^{-\frac{|c_B|^2}{s}} \l(\int_{B(x-c_B, r)} |\Omega(y)|^s \ddd y\r)^{\frac{1}{s}}.
\end{align*}
From $x\in B(c_B,r)$ we know that $0\in B(x-c_B,r)$ and then $B(x-c_B,r)\subset B(0,2r)$. Hence by using the polar coordinates we can deduce that
\begin{align*}
\l(\int_{B(c_B, r)} |\Omega(x-y)|^s e^{-|y|^2}\pi^{-\frac{d}{2}}\ddd y\r)^{\frac{1}{s}}
&\lesssim_a \pi^{-\frac{d}{2s}} e^{-\frac{|c_B|^2}{s}} \l(\int_{B(0, 2r)} |\Omega(y)|^s \ddd y\r)^{\frac{1}{s}} \\
&= \pi^{-\frac{d}{2s}} e^{-\frac{|c_B|^2}{s}} \l(\int_{0}^{2r}\rho^{d-1} \ddd\rho \int_{\bS^{d-1}} |\Omega(y')|^s \ddd \sigma(y')\r)^{\frac{1}{s}} \\
&= \pi^{-\frac{d}{2s}} e^{-\frac{|c_B|^2}{s}} \|\Omega\|_{L^{s}(\bS^{d-1})} \l(\frac{(2r)^d}{d}\r)^{1/{s}} \\
&\sim e^{-\frac{|c_B|^2}{s}} r^{\frac{d}{s}}\sim [\gamma(B)]^{\frac{1}{s}}.
\end{align*}
By taking this estimate into the inequality (\ref{DL_1998_prop1_1}) we obtain
\begin{align*}
\l|M_{\Omega, \beta}^a(f)(x)\r| &\lesssim \sup_{B\in\mathscr{B}_a(x)} \frac{1}{\gamma(B)^{1-\beta}}[\gamma(B)]^{1/{s}}\l(\int_B |f(y)|^{s'} \ddd\gamma(y)\r)^{1/{s'}} \\
&=\sup_{B\in\mathscr{B}_a(x)} \l(\frac{1}{\gamma(B)^{1-\beta s'}}\int_B |f(y)|^{s'} \ddd\gamma(y)\r)^{1/{s'}} \\
&=N_{\beta s', s'}^a f(x).
\end{align*}
From $1<s'<p<1/{\beta}, 1/q=1/p-\beta$ we know that $0<\beta s'<1, 1<\frac{p}{s'}<\frac{1}{\beta s'}$ and $\frac{1}{q/{s'}}=\frac{1}{p/{s'}}-\beta s'$. Then by Lemma \ref{DL_1998_lemma1} we finally prove the desired conclusion
\[\l(\int_{\bR^d} \l[M_{\Omega,\beta}^a f(x)\omega(x)\r]^q \ddd\gamma(x)\r)^{\frac{1}{q}}
\lesssim \l(\int_{\bR^d} \l[N_{\beta s', s'}^a f(x) \omega(x)\r]^q \ddd\gamma(x)\r)^{\frac{1}{q}}
\lesssim \l(\int_{\bR^d}|f(x)\omega(x)|^p \ddd\gamma(x)\r)^{\frac{1}{p}}.\]
\end{proof}

\begin{proof}[\textbf{Proof of Lemma \ref{LM_2020_lemma3.4}}]
The detailed proof is long but essentially the same as \cite[Lemma 3.4]{LM_2020}. For any $x\in\bR^d$, we divide the proof into two cases:

\textbf{Case (I)} $\gamma(B(x,am(x)))\leq \l(\frac{M_{\Omega, \beta_2}^{2a}(f)(x)}{M_{\Omega, \beta_1}^{2a}(f)(x)}\r)^{\frac{1}{\beta_2-\beta_1}}$.

By splitting the integral domain, together with the locally reverse doubling property of $\gamma$, we conclude that
\[|I_{\Omega, \beta}^{a}(f)(x)|\lesssim [\gamma(B(x,am(x)))]^{\beta-\beta_1}M_{\Omega, \beta_1}^{a}(f)(x)\leq [\gamma(B(x,am(x)))]^{\beta-\beta_1}M_{\Omega, \beta_1}^{2a}(f)(x).\]
then from the assumption in this case and $\beta>\beta_1$ we get the desired result.

\textbf{Case (II)} $\gamma(B(x,am(x)))> \l(\frac{M_{\Omega, \beta_2}^{2a}(f)(x)}{M_{\Omega, \beta_1}^{2a}(f)(x)}\r)^{\frac{1}{\beta_2-\beta_1}}$.

Actually we can use a similar argument in the previous case. By \cite[Lemma 2.11]{LSY_2014}, there exists an $r\in(0,a)$ such that
\[\frac{1}{2}\l(\frac{M_{\Omega,\beta_2}^{2a}(f)(x)}{M_{\Omega,\beta_1}^{2a}(f)(x)}\r)^{\frac{1}{\beta_2-\beta_1}} <\gamma(B(x,rm(x))) <\l(\frac{M_{\Omega,\beta_2}^{2a}(f)(x)}{M_{\Omega,\beta_1}^{2a}(f)(x)}\r)^{\frac{1}{\beta_2-\beta_1}}.\]
Then by dividing $I_{\Omega,\beta}^{a}(f)(x)$ into two parts as follows
\[I_1:=\int_{|x-y|<rm(x)}\frac{|f(y)|}{[\gamma(B(x,|x-y|))]^{1-\beta}}\ddd\gamma(y), \quad I_2:=\int_{rm(x)\leq|x-y|<am(x)}\frac{|f(y)|}{[\gamma(B(x,|x-y|))]^{1-\beta}}\ddd\gamma(y),\]
splitting the integral domain and using the locally reverse doubling property of $\gamma$ again, we can obtain
\[I_1\lesssim [\gamma(B(x,rm(x)))]^{\beta-\beta_1}M_{\Omega,\beta_1}^{2a}(f)(x), \quad I_2\lesssim [\gamma(B(x,rm(x)))]^{\beta-\beta_2}M_{\Omega,\beta_2}^{2a}(f)(x).\]
Combining these two cases, we deduce the desired result and complete the proof.
\end{proof}

\begin{proof}[\textbf{Proof of Theorem \ref{DL_1998_thm1}}]
Under the conditions of this theorem, by Lemma \ref{DL_1998_lemma3}, there exists an $\varepsilon>0$ such that
\[0<\varepsilon<\beta<\beta+\varepsilon<1,\quad 1/p>\beta+\varepsilon,\quad \omega^{s'}\in A(p/{s'},q_{\varepsilon}/{s'}, a), \quad \omega^{s'}\in A(p/{s'},\tilde{q}_{\varepsilon}/{s'},a),\]
where $1/{q_{\varepsilon}}=1/p-\beta-\varepsilon$ and $1/{\tilde{q}_{\varepsilon}}=1/p-\beta+\varepsilon$. Now set $l_1=2q_{\varepsilon}/q, l_2=2\tilde{q}_{\varepsilon}/q$, then $1/{l_1}+1/{l_2}=1$. For the $\varepsilon>0$ given above, by Corollary \ref{DL_1998_lemma4} and $\holder$'s inequality, we conclude that
\begin{align*}
\|I_{\Omega,\beta}^{a} f\|_{L^q(\bR^d, \gamma, \omega^q)}
&\lesssim \l(\int_{\bR^d} \l[M_{\Omega,\beta+\varepsilon}^{2a} (f)(x)\omega(x)\r]^{q/2}\l[M_{\Omega,\beta-\varepsilon}^{2a} (f)(x)\omega(x)\r]^{q/2} \ddd\gamma(x)\r)^{1/q} \\
&\leq \l(\int_{\bR^d} \l[M_{\Omega,\beta+\varepsilon}^{2a} (f)(x)\omega(x)\r]^{\frac{ql_1}{2}} \ddd\gamma(x)\r)^{\frac{1}{ql_1}} \l(\int_{\bR^d}\l[M_{\Omega,\beta-\varepsilon}^{2a} (f)(x)\omega(x)\r]^{\frac{ql_2}{2}} \ddd\gamma(x)\r)^{\frac{1}{ql_2}} \\
&=\l(\int_{\bR^d} \l[M_{\Omega,\beta+\varepsilon}^{2a} (f)(x)\omega(x)\r]^{q_{\varepsilon}} \ddd\gamma(x)\r)^{\frac{1}{2q_{\varepsilon}}} \l(\int_{\bR^d}\l[M_{\Omega,\beta-\varepsilon}^{2a} (f)(x)\omega(x)\r]^{\tilde{q}_{\varepsilon}} \ddd\gamma(x)\r)^{\frac{1}{2\tilde{q}_{\varepsilon}}}.
\end{align*}
Hence from Lemma \ref{DL_1998_lemma3}, Proposition \ref{DL_1998_prop1} and \cite[Corollary 2.9]{LM_2020}, we have
\[\l(\int_{\bR^d} \l[M_{\Omega,\beta+\varepsilon}^{2a} (f)(x)\omega(x)\r]^{q_{\varepsilon}} \ddd\gamma(x)\r)^{1/{2q_{\varepsilon}}} \lesssim \|f\|_{L^p(\bR^d, \gamma, \omega^p)}^{1/2}\]
and
\[\l(\int_{\bR^d}\l[M_{\Omega,\beta-\varepsilon}^{2a} (f)(x)\omega(x)\r]^{\tilde{q}_{\varepsilon}} \ddd\gamma(x)\r)^{1/{2\tilde{q}_{\varepsilon}}} \lesssim \|f\|_{L^p(\bR^d, \gamma, \omega^p)}^{1/2}.\]
Thereby we obtain the following desired result
\[\|I_{\Omega,\beta}^{a} f\|_{L^q(\bR^d, \gamma, \omega^q)}\lesssim \|f\|_{L^p(\bR^d, \gamma, \omega^p)}.\]
\end{proof}

\bigskip\bigskip

\subsection*{Acknowledgements} This work was supported by National Natural Science Foundation of China (Grant Nos. 11871452 and 12071473), Beijing Information Science and Technology University Foundation (Grant Nos. 2025031).

\bigskip\bigskip

\vspace{-0.4cm}{\footnotesize

\begin{flushleft}
	
		\vspace{0.3cm}\textsc{Boning Di\\School of Mathematical Sciences
		\\University of Chinese Academy of Sciences\\Beijing, 100049\\P. R. China}
	
	\vspace{0.3cm}\textsc{Qianjun He\\School of Applied
		Science\\Beijing Information Science and Technology University\\Beijing, 100192\\P. R. China}
	
	\emph{E-mail address}: \textsf{heqianjun16@mails.ucas.ac.cn}
	
	\vspace{0.3cm}\textsc{Dunyan Yan\\School of Mathematical Sciences
		\\University of Chinese Academy of Sciences\\Beijing, 100049\\P. R. China}
	
\end{flushleft}

\end{document}